\documentclass[reqno,12pt]{preprint}
\usepackage[full]{textcomp}
\usepackage{times}
\usepackage[cal=boondoxo]{mathalfa}
\usepackage{microtype}

\usepackage{comment}

\usepackage{amsmath}
\usepackage{amssymb}
\usepackage{amsfonts}
\usepackage{hyperref}
\usepackage{breakurl}

\usepackage{mathrsfs}
\usepackage{enumitem}

\usepackage{tikz}
\usepackage{dsfont}
\usepackage{mhequ} 
\usepackage{mhsymb} 
\usepackage{mhenvs} 
\usepackage{array}

\usepackage{wasysym}
\usepackage{centernot}

\setlist{noitemsep,nolistsep,leftmargin=1.7em}
\usetikzlibrary{snakes}
\usetikzlibrary{decorations}
\usetikzlibrary{positioning}
\usetikzlibrary{shapes}

\DeclareFontFamily{U}{mathx}{\hyphenchar\font45}
\DeclareFontShape{U}{mathx}{m}{n}{
      <5> <6> <7> <8> <9> <10>
      <10.95> <12> <14.4> <17.28> <20.74> <24.88>
      mathx10
      }{}
\DeclareSymbolFont{mathx}{U}{mathx}{m}{n}
\DeclareMathSymbol{\bigtimes}{1}{mathx}{"91}

\def\emptyset{{\centernot\ocircle}}

\definecolor{darkred}{rgb}{0.7,0.1,0.1}
\definecolor{darkblue}{rgb}{0.1,0.1,0.8}
\definecolor{darkgreen}{rgb}{0.1,0.7,0.1}
\def\CJ{\mathcal{J}}
\def\restr{\mathord{\upharpoonright}}

\providecommand{\figures}{false}
{ \ifthenelse{\equal{\figures}{false}} {#1}{\[ {\rm Figure \ missing !} \]} }{}
\def\id{\mathrm{id}}

\definecolor{connection}{rgb}{0.7,0.1,0.1}
\colorlet{symbols}{black!50}

\tikzset{
root/.style={circle,fill=black!50,inner sep=0pt, minimum size=3mm},
        dot/.style={circle,fill=black,inner sep=0pt, minimum size=1.5mm},
        dotred/.style={circle,fill=black!50,inner sep=0pt, minimum size=2mm},
        var/.style={circle,fill=black!10,draw=black,inner sep=0pt, minimum size=3mm},
        kernel/.style={semithick,shorten >=2pt,shorten <=2pt},
        kernels/.style={snake=zigzag,shorten >=2pt,shorten <=2pt,segment amplitude=1pt,segment length=4pt,line before snake=2pt,line after snake=5pt,},
        rho/.style={densely dashed,semithick,shorten >=2pt,shorten <=2pt},
           testfcn/.style={dotted,semithick,shorten >=2pt,shorten <=2pt},
        renorm/.style={shape=circle,fill=white,inner sep=1pt},
        labl/.style={shape=rectangle,fill=white,inner sep=1pt},
        xic/.style={very thin,circle,fill=symbols,draw=black,inner sep=0pt,minimum size=1.2mm},
        xi/.style={very thin,circle,fill=blue!10,draw=black,inner sep=0pt,minimum size=1.2mm},
        xix/.style={crosscircle,fill=blue!10,draw=black,inner sep=0pt,minimum size=1.2mm},
	xib/.style={very thin,circle,fill=blue!10,draw=black,inner sep=0pt,minimum size=1.6mm},
	xie/.style={very thin,circle,fill=green!50!black,draw=black,inner sep=0pt,minimum size=1.6mm},
	xid/.style={very thin,circle,fill=symbols,draw=black,inner sep=0pt,minimum size=1.6mm},
	xibx/.style={crosscircle,fill=blue!10,draw=black,inner sep=0pt,minimum size=1.6mm},
	kernels2/.style={very thick,draw=connection,segment length=12pt},
	not/.style={thin,circle,fill=symbols,draw=connection,fill=connection,inner sep=0pt,minimum size=0.5mm},
	>=stealth,
        }

\makeatletter
\def\DeclareSymbol#1#2#3{\expandafter\gdef\csname MH@symb@#1\endcsname{\tikz[baseline=#2,scale=0.15,draw=symbols,line join=round]{#3}}\expandafter\gdef\csname MH@symb@#1s\endcsname{\scalebox{0.7}{\tikz[baseline=#2,scale=0.15,draw=symbols,line join=round]{#3}}}}
\def\<#1>{\csname MH@symb@#1\endcsname}
\makeatother

\DeclareSymbol{Xi22}{0.5}{\draw (0,0) node[xi] {} -- (-1,1) node[not] {} -- (0,2) node[xi] {};}

\DeclareSymbol{Xi2}{-2}{\draw (0,-0.25) node[xi] {} -- (-1,1) node[xi] {};}
\DeclareSymbol{Xi3}{0}{\draw (0,0) node[xi] {} -- (-1,1) node[xi] {} -- (0,2) node[xi] {};}
\DeclareSymbol{Xi4}{2}{\draw (0,0) node[not] {} -- (-1,1) node[xi] {} -- (0,2) node[xi] {} -- (-1,3) node[xi] {};}
\DeclareSymbol{Xi2X}{-2}{\draw (0,-0.25) node[xi] {} -- (-1,1) node[xix] {};}
\DeclareSymbol{XXi2}{-2}{\draw (0,-0.25) node[xix] {} -- (-1,1) node[xi] {};}

\DeclareSymbol{IXi^2}{-1}{\draw (-1,1) node[xi] {} -- (0,0) node[not] {} -- (1,1) node[xi] {};}

\DeclareSymbol{XiX}{-2.8}{\node[xibx] {};}
\DeclareSymbol{tauX}{-2.8}{ \draw[kernels2] (0,0) node[xibx] {};}
\DeclareSymbol{Xi}{-2.8}{\node[xib] {};}
\DeclareSymbol{XiY}{-2.8}{\node[xie] {};}
\DeclareSymbol{XiZ}{-2.8}{\node[xid] {};}

\DeclareSymbol{I1Xi4}{2}{\draw (0,0) node[not] {} -- (-1,1) node[not] {}; \draw[kernels2] (-1,1) -- (0,2) ; 
\draw[kernels2] (-1,1) -- (-2,2) node[xi] {} ;
\draw (0,2) node[xi] {} -- (-1,3) node[xi] {};}

\DeclareSymbol{IXi}{0}{\draw (0,-0.25) node[not] {} -- (0,1.5) node[xi] {};}
\DeclareSymbol{IXiX}{0}{\draw (0,-0.25) node[not] {} -- (0,1.5) node[xix] {};}
\DeclareSymbol{IIXiX}{0}{\draw[kernels2] (-1,1) -- (0,2) node[xi] {}; \draw (0,-0.25) node[not] {} -- (-1,1) node[not] {}; }

\DeclareSymbol{IXi2}{0}{\draw (0,-0.25) node[not] {} -- (-1,1) node[xi] {} -- (0,2) node[xi] {};}
\DeclareSymbol{IIXi^2}{-4}{\draw (0,-1.5) node[not] {} -- (0,0);
\draw[kernels2] (-1,1) node[xi] {} -- (0,0) node[not] {} -- (1,1) node[xi] {};}
\DeclareSymbol{I'IXi^2}{-4}{\draw[kernels2] (0,-1.5) node[not] {} -- (0,0);
\draw[kernels2] (-1,1) node[xi] {} -- (0,0) node[not] {} -- (1,1) node[xi] {};}

\DeclareSymbol{I1Xi4bc}{2}{\draw (0,0) node[not] {} -- (-1,1) node[not] {}; \draw[kernels2] (-1,1) -- (0,2) ; 
\draw[kernels2] (-1,1) -- (-2,2) node[xi] {} ; \draw[kernels2] (0,2) node[not] {} -- (-1,3) node[xi] {};\draw[kernels2] (0,2)  -- (1,3) node[xi] {};
}

\DeclareSymbol{Xi4ca}{0}{\draw (0,1) -- (-1.3,2.2) node[xi] {};\draw (0,-0.25) node[not] {} -- (0,1) ; \draw[kernels2] (0,1) node[not] {} -- (1.3,2.2) node[xi] {};
\draw[kernels2] (0,1) {} -- (0,2.7) node[xi] {};
}

\DeclareSymbol{Xi4c}{0}{\draw (0,1) -- (0.8,2.2) node[xi] {};\draw (0,-0.5) node[not] {} -- (0,1) node[xi] {} -- (-0.8,2.2) node[xi] {};}

\DeclareSymbol{I1Xi4b}{2}{\draw (0,0) node[not] {} -- (-1,1) node[xi] {} -- (0,2) ; \draw[kernels2] (0,2) node[not] {} -- (-1,3) node[xi] {};\draw[kernels2] (0,2)  -- (1,3) node[xi] {};
}

\def\CP{\mathcal{P}}
\def\CW{\mathcal{W}}
\def\CG{\mathcal{G}}
\def\CJ{\mathcal{J}}
\def\CA{\mathcal{A}}

\def\CC{\mathcal{C}}

\def\CM{\mathcal{M}}
\def\CT{\mathcal{T}}
\def\RR{\mathfrak{R}}

\def\Labe{\mathfrak{e}}
\def\m{\mathfrak{m}}
\def\Labn{\mathfrak{n}}

\def\Deltam{\Delta^{\!-}}

\def\Deltap{\Delta^{\!+}}

\def\scal#1{\langle#1\rangle}

\def\Trees{\mathfrak{T}}

\def\CS{\mathcal{S}}
\def\CR{\mathcal{R}}
\def\CF{\mathcal{F}}

\newenvironment{DIFnomarkup}{}{} 

\newfont{\indic}{bbmss12}

\def\ex{\mathrm{ex}}

\def\PPi{\boldsymbol{\Pi}}

\def\one{\mathbf{1}}
\def\eps{\varepsilon}

\def\TT{\mathscr{T}}
\def\MM{\mathscr{M}}

\def\Deltap{\Delta^{\!+}}
\def\Deltapm{\Delta^{\!\pm}}
\def\Deltam{\Delta^{\!-}}

\def\PPi{\boldsymbol{\Pi}}

\def\id{\mathrm{id}}

\def\RR{\mathfrak{R}}
\def\br{\penalty0}

\def\simnot{\stackrel{\vbox to 0.15em{\hbox{\kern0.07em$^\circ$}}}{\sim}}

\begin{document}

\title{The motion of a random string}
\author{M. Hairer}
\institute{Mathematics Institute, The University of Warwick}

\maketitle

\begin{abstract}
We review a series of forthcoming results leading to the construction of 
a natural evolution on the space of loops with values in
a Riemannian manifold. In particular, this clarifies the algebraic structure
of the renormalisation procedures appearing in the context of the
theory of regularity structures.
\end{abstract}

\section{Introduction}

This work is partially 
motivated by Funaki's attempt \cite{FunakiOld,Funaki} to construct a natural evolution
on the space of loops with values in a manifold.
Given a compact smooth Riemannian manifold $\CM$ with metric $g$, 
write $\CL^\infty \CM$ for the space of all smooth loops in $\CM$, i.e.\
simply the space of all smooth functions $u \colon S^1 \to \CM$.
We also write $\CL\CM$ for the space of all loops, where we only impose that $u$ is
continuous.
The energy of a 
loop $u\in \CL^\infty\CM$ is given by
\begin{equ}[e:energy]
E(u) = {1\over 2}\int_{S^1} g_{u(x)}(\d_x u(x),\d_x u(x))\,dx\;.
\end{equ}
The aim of this note is to discuss the construction of  the natural
Langevin dynamic associated to $E$. In other words, one would like to 
build a Markov process $u$ with values in $\CL\CM$ with the
following two properties:
\begin{claim}
\item[1.] The measure on loops formally given by $\exp(-2E(u))\,Du$ is 
invariant for the process $u$.  
\item[2.] The evolution is local in the sense that, when written as an 
evolution equation, all of the terms appearing in its right hand side
are purely local.
\end{claim}
The meaning of the first property is not clear a priori since there is of course
no ``Lebesgue measure $Du$ on loop space'', but a natural way of interpreting
it is as the Brownian bridge measure on $\CM$. See \cite{MR780661,MR1698956}
for proofs that natural approximations of $\exp(-2E(u))\,Du$ do indeed
converge to the Brownian bridge measure, or to close relatives thereof.
The second property is natural in view of the fact that 
the integrand appearing in \eqref{e:energy} depends on $u$ in a  local
way. 

\subsection{Derivation of the main equation}

Recall that in the finite-dimensional case, given 
a smooth potential $V \colon \R^d \to \R$ with sufficiently ``nice'' behaviour
at infinity, the Langevin dynamic associated to $V$ is given by
the stochastic differential equation
\begin{equ}[e:Langevin]
\dot x = -\nabla V(x) + \xi\;,
\end{equ} 
where $\xi$ denotes white noise, i.e.\ the generalised $\R^d$-valued stochastic
process with 
\begin{equ}[e:defXi]
\E \xi_i(s)\xi_j(t) = \delta_{ij}\delta(t-s)\;.
\end{equ}
In order to mimic this, we therefore first need a natural
notion of gradient for $E$. 
For this, note that the metric $g$ induces a Riemannian structure $\hat g$ on 
$\CL^\infty \CM$ in the following way. An element $h \in T_u \CL^\infty\CM$
of the tangent space of $\CL^\infty \CM$ at $u \in \CL^\infty \CM$
is given by a function $h \colon S^1 \to T\CM$ is such that $h(x) \in T_{u(x)}\CM$ for
every $x \in S^1$. One then defines a ``Riemannian metric'' $\hat g$ on 
$T \CL^\infty\CM$ by setting 
\begin{equ}
\hat g_u(h,h) = \int_{S^1} g_{u(x)}(h(x),h(x))\,dx\;,
\end{equ}
for every $h \in T_u \CL^\infty\CM$. As usual, one can then define the gradient 
$\nabla E(u)$ as the unique element of $T_u\CL^\infty \CM$ such that 
\begin{equ}
\hat g_u(h, \nabla E(u)) = dE(h)\;,
\end{equ}
for every $h \in T_u\CL^\infty \CM$, where $dE_u(h)$ is the differential of $E$
in the direction $h$ evaluated at $u$. This yields a particular case of the Eells-Sampson
Laplacian \cite{ES}.

In view of \eqref{e:Langevin}, it is then natural to try to define
the Langevin equation associated to \eqref{e:energy} as
\begin{equ}[e:LangevinLoop]
\dot u = -\nabla E(u) + \xi\;,
\end{equ}
where $\xi$ is a suitable ``white noise on loop space''. In particular, in the
simplest possible case where $\CM = \R$ (with its standard metric), one
obtains the stochastic heat equation $\dot u = \d_x^2 u + \xi$, for $\xi$ 
a space-time white noise, the invariant measure of which is indeed given by the
Brownian bridge. 
In general however, it is much less clear what the symbol $\xi$ appearing
in \eqref{e:LangevinLoop}
latter actually means since at any fixed time $\xi$ should belong to $T_u \CL\CM$
(or at least a suitable generalisation thereof), which itself
varies with time.

Note that a characterisation of the noise $\xi$ appearing in \eqref{e:Langevin} is that it
is the (generalised) Gaussian process with Cameron-Martin 
space given by $L^2(\R_+, \R^d)$ endowed with its usual scalar product. This hints
at a way of adapting this definition to the situation at hand.
Indeed, for any given smooth function $u \colon \R_+ \to \CL^\infty \CM$ 
(which of course the process $u$ we are interested in really isn't!), 
we obtain a scalar product on the space $L^2_u$ of functions
$h \colon \R_+ \to T \CL^\infty\CM$ with $h(t) \in T_{u(t)} \CL\CM$ by
\begin{equ}[e:defScalar]
\scal{h,\bar h}_u = \int_\R \int_{S^1}\scal{h(t,x),\bar h(t,x)}_{u(t,x)}\,dx\,dt\;.
\end{equ}
then a natural definition for $\xi$ mimicking \eqref{e:defXi} is that it is the
(generalised) Gaussian process with Cameron-Martin space given by $L^2_u$ endowed with the
scalar product \eqref{e:defScalar}. It turns out that this process can 
be constructed more explicitly as follows. 
Take a finite collection $\{\sigma_i\}_{i=1}^m$ of vector
fields on $\CM$ with the property that, for every $u \in \CM$ and $h, \bar h \in T_u \CM$,
one has the identity
\begin{equ}
g_u(h,\bar h) = \sum_{i=1}^m g_u(h,\sigma_i(u))g_u(\bar h,\sigma_i(u))\;,
\end{equ}
and let $\{\xi_i\}_{i=1}^m$ be a collection of independent space-time white noises.
Then, one can verify that the Gaussian process $\xi$ given by $\xi(t,x) = \sum_{i=1}^m 
\sigma_i(u(t,x))\,\xi_i(t,x)$ does indeed have as its Cameron-Martin
space $L^2_u$ with scalar product \eqref{e:defScalar}.

This suggests that \eqref{e:LangevinLoop} should really be interpreted as
\begin{equ}
\dot u = -\nabla E(u) + \sum_{i=1}^m\sigma_i(u)\,\xi_i\;.
\end{equ}
In local coordinates, this can be written as
\begin{equ}[e:main]
\dot u^\alpha = \d_x^2 u^\alpha + \Gamma^\alpha_{\beta\gamma}(u)\, \d_x u^\beta\d_x u^\gamma + \sigma_i^\alpha(u)\,\xi_i\;,
\end{equ}
where Einstein's convention of summation over repeated indices is implied
and $\Gamma^\alpha_{\beta\gamma}$ are the Christoffel symbols for the
Levi-Civita connection of $(\CM,g)$. We thus derived some kind of multi-component
version of the KPZ equation, similar to the one studied in \cite{FSS,Spohn}, 
but with both $\Gamma$ and $\sigma$
depending on the solution $u$ itself.
This also shows where the problem lies.
Since we expect the stationary measure of \eqref{e:main} to be closely related to
the Brownian bridge, we expect typical solutions to be not even 
${1\over 2}$-H\"older continuous as a function of the $x$ coordinate. As a consequence,
the term $\d_x u^\beta\d_x u^\gamma$ involves the product of two distributions, so that
it is not clear \textit{a priori} how to define it. Even if we were to somehow define it,
it still wouldn't be clear how to multiply the resulting distribution with the
very irregular function $\Gamma^\alpha_{\beta\gamma}(u)$. 

Similarly, while one would
expect to be able to define the terms $\sigma_i^\alpha(u)\,\xi_i$ via It\^o integration, 
this would certainly lead to problems since It\^o integration famously does not preserve
the chain rule, which makes it ill-suited for problems on manifolds. On the other hand,
there appears to be no Stratonovich integration in this context because the 
corresponding It\^o correction term is infinite, see for example \cite{WongZakai}.

\subsection{Main result}

In the present work, we take a somewhat pragmatic approach by postulating
that we call a process a ``solution'' to \eqref{e:main} if it is the limit
of a natural approximation procedure. More precisely, we fix a function
$\rho\in \CC_0^\infty(\R^2)$ integrating to $1$ and we set
\begin{equ}[e:scalerho]
\xi_i^{(\eps)} = \rho^{(\eps)}\star \xi \;,\qquad \rho^{(\eps)}(t,x) = \eps^{-3} \rho(t/\eps^2, x/\eps)\;.
\end{equ}
It is then natural to ask whether solutions to \eqref{e:main} with $\xi$
replaced by $\xi^{(\eps)}$ admit limits as $\eps \to 0$.
Of course, even if this turns out to be the case, one might lose uniqueness,
but as we will see shortly this can be recovered to some extent. In order to
formulate our result, we impose an additional constraint on the vector fields
$\sigma_i$, namely we assume that 
\begin{equ}
\sum_{i=1}^m \nabla_{\sigma_i} \sigma_i = 0\;,
\end{equ}
where $\nabla$ denotes the covariant derivative (with respect to the Levi-Civita
connection). An equivalent way of formulating this is that, viewing
the $\sigma_i$ as first-order differential operators, $\sum_i \sigma_i^2$ is
equal to the Laplace-Beltrami operator. This should be viewed as a ``centering''
condition, which can always be satisfied. For example, it suffices to consider an
isometric embedding of $\CM$ as a submanifold of $\R^m$ and to then take for $\sigma_i(u)$ 
the orthogonal projection of the $i$th canonical basis vector of $\R^m$ onto
$T_u\CM$.

The main result on which we report here is a consequence of the more general
results exposed in \cite{YvainLorenzo,Ajay}, relying on \cite{Regularity,HQ}, and
can be formulated as follows.

\begin{theorem}\label{theo:main}
Let $u^{(\eps)}$ be the solution to \eqref{e:main} with $\xi_i$ replaced by $\xi_i^{(\eps)}$
and with some fixed H\"older continuous initial condition $u_0$.
Then, as $\eps \to 0$, $u^{(\eps)}$ converges in probability to a Markov process $u$.
Furthermore, although $u$ might depend on the choice of $\rho$, the set of possible
limits can be parametrised smoothly by finitely many parameters.
\end{theorem}

\begin{remark}
If $\CM$ is a symmetric space, then the limit can be shown to be independent of $\rho$. 
In general, we believe that the set of possible limits for approximations of the
type discussed here is actually one-dimensional, but we do not have a proof of this.
We also believe that the law of the limit is independent of the choice of vector fields
$\sigma_i$, but again we have no proof of this.
\end{remark}

\begin{remark}
A natural process on loop space leaving invariant the Brownian bridge measure on
a manifold was constructed by a number of authors in the nineties,
see for example \cite{Driver1,DriverRockner,Driver2,Stroock,Norris}. This process 
is different since, in the particular case of $\CM = \R^d$, it would
correspond to the Ornstein-Uhlenbeck process from Malliavin calculus, rather than
the stochastic heat equation. In particular, it is not ``local'' in the sense that 
its driving noise has non-trivial spatial correlations.
\end{remark}

To some extent, we can probably learn more from the proof of Theorem~\ref{theo:main} 
than from the result itself, although a more systematic study 
of the properties of the limiting process,
especially in the small noise limit, could be very interesting, especially in view of 
\cite{Witten,Bismut}. 
We will assume in the sequel that the reader has at least some familiarity with the 
theory of regularity structures as
developed in \cite{Regularity} and surveyed for example in \cite{Review,Review2}.
The proof follows the same methodology as developed in
\cite{Regularity,HQ,WongZakai}, which can be summarised as follows:
\begin{claim}
\item[1.] Build a regularity structure $\TT = (\CT,\CG_+)$ that is sufficiently rich to be able
to formulate \eqref{e:main} as a fixed point problem in some space $\CD^{\gamma,\eta}$
of sufficiently regular modelled distributions and build a local solution theory.
\item[2.] Exhibit a sufficiently large subgroup $\CG_-$ of the ``renormalisation group'' $\RR$
associated to $\TT$ as in \cite[Def.~8.41]{Regularity}.
\item[3.] Find a sequence of elements $M_\eps \in \CG_-$ such that the sequence of admissible 
models $\hat \Pi^{(\eps)}$ obtained by taking the canonical lift $\Pi^{(\eps)}$ of $\{\xi_i^{(\eps)}\}_{i=1}^m$ 
and acting on it with $M_\eps$ converges to a limiting model $\hat \Pi$.
\item[4.] Show that $M_\eps$ can be chosen in such a way that the solution to the fixed point
problem constructed in 1.\ with model $\hat \Pi^{(\eps)}$ is the same (once higher-order information 
is discarded) as that of \eqref{e:main} with $\xi_i$ replaced by $\xi_i^{(\eps)}$.
\end{claim}
It follows immediately from the main results of \cite{Regularity} that if these four steps
can be carried out successfully, then Theorem~\ref{theo:main} follows.

Let us now review how these steps can be performed in the context of the problem at hand.
Step 1.\ in this list is generic and was performed in \cite{Regularity} for a large class of stochastic PDEs which
in particular include \eqref{e:main}, so we will not dwell on it. Step 2.\ is purely
algebraic. The problem is that while elements of $\RR$ are described by linear maps on $\CT$
and one can easily make an ``educated guess'' of a group $\CG_-$ of such linear maps that 
is relevant later on, checking whether a given map belongs to $\RR$ using the
characterisation given in \cite{Regularity} is tedious. In the
case of \cite{Regularity,HQ,WongZakai}, the relevant groups were of very low dimension, 
so that this verification could easily be done by brute force. In the present situation, even after
taking simplifications arising from symmetries into account, 
the dimension of $\CG_-$ is still $71$, so that such a brute
force verification would only be practical with the help of a computer. In \cite{YvainLorenzo}, 
we therefore obtain a general result showing that those linear maps arising from the
``educated guesses'' alluded to earlier do indeed always belong to $\RR$, and a presentation of this
result is the content of Section~\ref{sec:algebra}.
Step 3.\ on the other hand is purely analytical. Again, in the abovementioned previous works, 
this step was reduced to moment bounds for a finite number (typically about $5$ or $6$) 
of stochastic processes, and these were then obtained separately for each of these processes.
Again, this could in principle be done here but would be impractical since the number of processes
that would require moment bounds in this case is well over one hundred. 
In the article \cite{Ajay}, we therefore obtain a general estimate that yields suitable
moment bounds for a very large class of processes of this type, thus bypassing the need to 
consider these processes separately. This bound is presented in Section~\ref{sec:analysis}.
At this stage, without any additional information but exploiting the $x \leftrightarrow -x$ symmetry
of the equation, one deduces from the general theory of \cite{Regularity} that there exist \textit{finitely many} 
$\eps$-dependent constants $c_j^{(\eps)}$ and $\eps$-independent vector fields $h_j$ such that 
the solutions to
\begin{equ}
\dot u^\alpha = \d_x^2 u^\alpha + \Gamma^\alpha_{\beta\gamma}(u)\, \d_x u^\beta\d_x u^\gamma + c_j^{(\eps)} h_j^\alpha(u) + \sigma_i^\alpha(u)\,\xi_i^{(\eps)}\;,
\end{equ}
converge to a limit $u$ that is furthermore independent of the mollifier $\rho$.
In order to complete the proof of Theorem~\ref{theo:main}, it therefore remains to show that 
for every $j$, either $c_j^{(\eps)}$ converges to a finite limit as $\eps \to 0$ or
$h_j = 0$. This result is the content of forthcoming work and 
will not be reviewed presently.

\section{Algebraic aspects}
\label{sec:algebra}

Recall that the regularity structure canonically associated to \eqref{e:main}
comprises symbols $\{\Xi_i\}_{i=1}^m$, $\{X^k\}_{k \in \N^2}$ representing the
driving noises and Taylor monomials respectively, as well as abstract 
integration maps $\{\CI_k\}_{k \in \N^2}$ representing the heat kernel and its derivatives, 
and an associative and commutative product. We write $\CF$ for the free vector space generated by 
all formal expressions built from these operations, modulo the usual identifications, namely that 
$\one \eqdef X^0$ is neutral for the product, that $X^k X^\ell = X^{k+\ell}$, 
and that $\CI_\ell(X^k) = 0$.
Writing $|k| = 2k_0 + k_1$ for $k \in \N^2$, we assign real-valued degrees to these objects by
\begin{equ}[e:grading]
|\Xi_i| = -{3\over 2} - \kappa\;,\qquad
|X^k| = |k|\;,\qquad
|\CI_\ell| = 2 - |\ell|\;,
\end{equ}
respectively, for $\kappa > 0$ a parameter with a sufficiently small value
($\kappa = 1/100$ will do). We will use the shorthand notations
$\CI' = \CI_{(0,1)}$, $X_0 = X^{(1,0)}$ and $X_1 = X^{(0,1)}$.
The fact that the time coordinate $X_0$ has
degree $2$ rather than $1$ reflects the fact that we endow space-time with the
parabolic scaling, which is consistent with the scaling used in \eqref{e:scalerho}.

Of course, not all the formal expressions built from these symbols and operations 
are useful for our problem. As in \cite{Regularity}, we define collections $\CU$, $\CU'$
and $\CW$ of formal expressions as the smallest collections such that
$\{\Xi_i\}_{i=1}^m \cup \{X^k\}_{k \in \N^2} \subset \CW$, such that
\begin{equ}
\tau \in \CW \qquad \Rightarrow\qquad \CI(\tau) \in \CU\quad\&\quad \CI'(\tau) \in \CU'\;,
\end{equ}
and such that, for every $k \ge 0$,
every collection $\{\tau_i\}_{i=1}^k$ with $\tau_i \in \CU$
for every $i$, and every $\sigma_1,\sigma_2 \in \CU'$, one has
\begin{equ}
\{\tau,
\tau \Xi_j,
\tau \sigma_1,
\tau \sigma_1\sigma_2\} \subset \CW\;,\qquad \tau = \tau_1\cdots\tau_k\;.
\end{equ}
Note that by construction one has $\CU \cup \CU' \subset \CW$.

The space $\CT$ is then defined as the free vector space generated by 
$\CW$, endowed with the $\R$-grading determined by
\eqref{e:grading}. As long as $\kappa < {1\over 2}$, one can verify that,
for every $\gamma \in \R$, the subspace of $\CT$ generated by elements of degree at
most $\gamma$ is finite-dimensional. 
In order to describe $\CG_+$, we consider the free vector space $\hat\CT_+$ generated by all formal 
expressions of the type
$X^\ell \prod_{i=1}^N \CI_{k_i}(\tau_i)$ for some $N \ge 0$, some $\ell, k_i \in \N^2$ and some
$\tau_i \in \CW$. Note that 
$\hat \CT_+$ is an algebra, which is not the case for $\CT$ since the product of two expressions
in $\CW$ does not necessarily belong to $\CW$. We then write $\CT_+$ for the quotient of $\hat \CT_+$ by the
ideal $\CJ_{-} \subset \hat \CT_+$ 
generated by all $\CI_{k}(\tau)$ with $|\CI_{k}(\tau)| \le 0$.
 
Consider now the linear map $\Deltap \colon \CT \to \CT \otimes \CT_+$ given by 
\begin{equ}[e:deltap1]
\Deltap X_i = X_i \otimes \one + \one \otimes X_i\;,\qquad
\Deltap \Xi_i = \Xi_i \otimes \one\;,
\end{equ}
and then recursively by
\begin{equ}[e:recDelta]
\Deltap \CI_k(\tau) = (\CI_k \otimes \id)\Deltap \tau
+ \sum_{\ell \in \N^2} {X^\ell \over \ell!} \otimes \CI_{k+\ell}(\tau)\;,\qquad
\Deltap (\tau \bar \tau) = \Deltap \tau\,\Deltap \bar \tau\;.
\end{equ}
The sum appearing in this expression is actually finite since all but finitely many
of the summands are zero in $\CT \otimes \CT_+$. Note also that 
we can define $\Deltap \colon \CT_+ \to \CT_+ \otimes \CT_+$ in formally exactly the same way
(and this is compatible with the quotienting procedure used to produce $\CT_+$), except that in
that case the maps $\CI_k$ appearing in \eqref{e:recDelta} are interpreted as linear maps
from $\CT$ to $\CT_+$. The following was shown in \cite[Sec.~8]{Regularity}.

\begin{proposition}\label{prop:Hopf+}
There exists an algebra morphism $\CA_+ \colon \CT_+ \to \CT_+$ so that 
$(\CT_+,\cdot,\br\Deltap,\CA_+)$ is a Hopf algebra, and 
the map $\Deltap$, viewed as a map $\CT \to \CT \otimes \CT_+$, turns
$\CT$ into a right comodule for $\CT_+$.
\end{proposition}

The morphism $\CA_+$ is defined uniquely by the fact that $\CA_+ X_i = - X_i$ and
\begin{equ}[e:defA+]
\CA_+ \CI_k(\tau) = -\sum_{\ell\in \N^2}{(-X)^\ell \over \ell!} \CM_+ \bigl(\CI_{k+\ell}\otimes \CA_+\bigr)\Deltap \tau\;,
\end{equ}
where $\CM_+ \colon \CT_+\otimes \CT_+ \to \CT_+$ denotes the product.

Denote by $\CG_+$ the group of characters $g \colon \CT_+ \to \R$ endowed with
\begin{equ}[e:defProd]
f \circ g = (f\otimes g)\Deltap\;,\qquad 
f^{-1} = f\CA_+\;,
\end{equ}
as well as the identity element $e$ given by $e(\one) = 1$ and $e(\tau) = 0$ for symbols
$\tau \neq \one$.
The comodule structure of $\CT$ for $\CT_+$ then yields a natural 
action of $\CG_+$ onto $\CT$ from the left via
\begin{equ}
g \mapsto \Gamma_g \;,\qquad \Gamma_g = (\id \otimes g)\Deltap\;.
\end{equ}
Since furthermore $\Deltap$ preserves total degree and every basis vector in 
$\CT_+$ has strictly positive degree except for $\one$, this 
yields a regularity structure $\TT = (\CT,\CG_+)$.
As shown in \cite{YvainLorenzo} and although it may not look so at first sight, this construction is 
equivalent to the one given in \cite[Sec.~8]{Regularity} modulo a change of basis
for $\CT_+$.

Furthermore, given a collection $\{\xi_i^{(\eps)}\}_{i=1}^m$ of smooth functions, 
we have a canonical lift to a model $(\Pi^{(\eps)},f^{(\eps)})$ for $\TT$ in the following way.
First, we define a linear map $\PPi^{(\eps)}\colon \CT \to \CC^\infty$ by
\minilab{e:model}
\begin{equs}[2]
\PPi^{(\eps)} \Xi_i &= \xi_i^{(\eps)}\;,&\qquad
\PPi^{(\eps)} \tau \bar \tau 
&= 
(\PPi^{(\eps)} \tau)\cdot (\PPi^{(\eps)} \bar \tau)\;,\label{e:canonical}\\
\PPi^{(\eps)} X^k\tau &= x^k \PPi^{(\eps)}\tau\;,&\qquad
\PPi^{(\eps)} \CI_k \tau &=  D^k K\star \PPi^{(\eps)}\tau\;,\label{e:admissible}
\end{equs}
where the kernel $K$ is a truncated heat kernel as in \cite{Regularity}, $\star$
is convolution in $\R^2$, and
$x^k$ denotes the function $x \mapsto x^k$. In general, we define

\begin{definition}
A linear map $\PPi\colon \CT \to\CC^\infty$ is \textit{admissible} if
\eqref{e:admissible} holds.
\end{definition}

Given any admissible linear map $\PPi\colon \CT \to \CC^\infty$, there is a natural
way of assigning to it a collection of characters $g^+_z(\PPi)\colon \hat \CT_+ \to \R$
by setting $g^+_z(\PPi)X_i = z_i$, $g^+_z(\PPi)\CI_k(\tau) = (D^k K\star \PPi\tau)(z)$,
and then extending this multiplicatively.
Write now $\hat\CA_+\colon \CT_+ \to \hat \CT_+$ for the unique algebra morphism
such that $\hat\CA_+ X_i = - X_i$ and
\begin{equ}
\hat\CA_+ \CI_k(\tau) = -\sum_{\ell}{(-X)^\ell \over \ell!} \pi_+ \hat \CM_+ \bigl(\CI_{k+\ell}\otimes \hat\CA_+\bigr)\Deltap \tau\;,
\end{equ}
where this time the $\CI_k$ are interpreted as maps $\CT \to \hat \CT_+$ and 
$\pi_+$ is the projection in $\hat \CT_+$ onto the terms of positive homogeneity.
(Similarly to above, $\hat \CM_+$ is the product in $\hat\CT_+$.)
Comparing this to the definition of $\CA_+$ given in \eqref{e:defA+}, we see that
it is virtually identical, with the exception of the appearance of the projection $\pi_+$
and the fact that the two operators do not act on the same spaces. 

Given $\PPi$, we then define  $\Pi_z\colon \CT \to \CC^\infty$ and 
$f_z\in \CG_+$ by
\begin{equ}[e:defModel1]
f_z =  g_z^+(\PPi)\hat\CA_+\;,\qquad \Pi_z = \bigl(\PPi \otimes f_z\bigr)\Deltap \;,
\end{equ}
for every $z \in \R^2$. We also define $\Gamma_{\!z\bar z}\colon \CT \to \CT$ by
\begin{equ}[e:defModel2]
\Gamma_{\!z\bar z} = \bigl(\id\otimes f_z\CA \otimes f_{\bar z}\bigr)(\Deltap \otimes \id)\Deltap 
= (\gamma_{z\bar z} \otimes \id)\Deltap\;,
\end{equ}
with $\gamma_{z\bar z} = (f_z\CA \otimes f_{\bar z})\Deltap$,
so that $\Pi_z \Gamma_{\!z\bar z} = \Pi_{\bar z}$ and we henceforth write $\CZ$ for the map
\begin{equ}
\CZ \colon \PPi \mapsto (\Pi,\Gamma)\;,
\end{equ}
given by \eqref{e:defModel1} and \eqref{e:defModel2}.
It is then possible to verify that if $\PPi^{(\eps)}$ is given by
\eqref{e:model}, then $\CZ(\PPi^{(\eps)})$ is 
indeed an admissible model for $\TT$. For a generic admissible
linear map $\PPi\colon \CT \to \CC^\infty$
however this is not necessarily true.

Writing $\MM$ for the space of admissible models for $\TT$ that are
periodic in the spatial variable, one can build a 
solution map $\CS \colon \MM \times \CL^\alpha\CM \to \bar \CC(\R_+, \CL^\alpha \CM)$ 
(here one should take $\alpha \in (0,{1\over 2}-\kappa)$, $\CL^\alpha\CM$ denotes the
space of $\alpha$-H\"older continuous loops in $\CM$, and 
$\bar \CC(\R_+,X)$ denotes the space of continuous functions with values in the
metric space $X$, up 
to some explosion time at which they leave every bounded region of $X$)
with the following two properties:
\begin{claim}
\item If $\PPi^{(\eps)}$ is the canonical lift for 
some smooth functions $\xi_i^{(\eps)}$, then $\CS(\CZ(\PPi^{(\eps)}), u_0)$
is the maximal solution to \eqref{e:main} with initial condition $u_0$ and $\xi_i$
replaced by $\xi_i^\eps$.
\item The map $\CS$ is locally Lipschitz continuous in both of its arguments.
\end{claim} 
This shows that if it were the case that $\CZ(\PPi^{(\eps)})$ converges
to some limiting model in $\MM$ as $\eps \to 0$, then Theorem~\ref{theo:main}
would follow at once. Unfortunately, this is simply not the case. We would therefore
like to be able to ``tweak'' this model in such a way that it remains an admissible
model but has a chance of converging as $\eps \to 0$.

A natural way of ``tweaking'' $\PPi^{(\eps)}$ is to compose it with some
linear map $M\colon \CT \to \CT$.
This naturally leads to the following question: what are the linear maps
$M$ which are such that if $\CZ(\PPi)$ is an admissible
model, then $\CZ(\PPi M)$ is also an admissible model? More precisely, 
we give the following definition.

\begin{definition}
A linear map $M \colon \CT \to \CT$ is an \textit{admissible renormalisation procedure}
if, for every $\PPi \colon \CT \to \CC^\infty$ such that 
$\CZ(\PPi) \in \MM$, one has $\CZ(\PPi M) \in \MM$.
\end{definition}

In \cite[Sec.~8.3]{Regularity}, 
a rather
indirect characterisation of renormalisation procedures $M$ is given.
The aim of the remainder of this section is to give an explicit description of
a very large class of such $M$.
In order to describe these maps, we recall that the elements of $\CW$ can best 
be described by rooted trees with additional decorations in the following way. 
We have $m+1$ types of edges: $m$ of them represent the symbols $\Xi_i$, and they always
touch a leaf. The last type represents $\CI$ and is decorated with an $\N^2$-valued label $k$,
thus representing $\CI_k$. Nodes are also decorated with $\N^2$-valued labels, 
representing factors of
$X^k$. Finally, multiplication is represented by concatenation at the root.
If we draw \tikz[baseline=-2.5] \node[not] {}; for nodes, \<Xi> for $\Xi_i$ (this looks like a node, but can unambiguously be 
interpreted as an edge protruding from the center of the disk since we postulated that 
these edges always terminate in a leaf),  \tikz[baseline=-2.5] \draw[symbols] (0,0) -- (0.4,0);
for $\CI$ and \tikz[baseline=-2.5] \draw[kernels2] (0,0) -- (0.4,0); for $\CI'$, we do for example
have
\begin{equ}
\CI \bigl(\CI'(\Xi_i)\CI'(\Xi_j)\bigr) = \<IIXi^2>\;,\qquad
\CI \bigl(\CI'(\Xi_i)\CI'(\Xi_j\CI(\Xi_k))\bigr) = \<I1Xi4>\;,
\end{equ}
etc, with the understanding that \<Xi> denotes the relevant $\Xi_i$. 
The expression $\CI' \bigl(X^k\CI'(\Xi_i)\CI'(\Xi_j)\bigr)$ for example would be represented
by \<I'IXi^2>, but with a decoration $k$ at the center node where the three bold lines meet. 
This gives a bijection between canonical basis vectors of $\CF$ and
triples $(F,\Labn,\Labe)$ where $F= (V_F,E_F)$ is a rooted trees with edge types in $\{\CI,\Xi_1,\ldots,\Xi_m\}$ subject to the restrictions described above,
$\Labn\colon V_F \to \N^2$, and $\Labe\colon E_F \to \N^2$.

By analogy with the BPHZ renormalisation procedure \cite{BP,Hepp,Zimmermann}, one
would like to consider renormalisation maps that consist in ``contracting subtrees''.
In order to formalise such an operation, consider a tree $T = (V_T,E_T)$, as well as 
a subforest $A = (V_A,E_A) \subset T$, i.e.\ an arbitrary subgraph of $T$ which contains
no isolated vertices. We then write $\CR_A T$ for the 
tree obtained by contracting the connected components of $A$ in $T$.

We also write $\hat \CT_-$ for the free commutative algebra
generated by $\CW$ 
and $\CJ_+ \subset \hat \CT_-$ for the ideal 
generated by $\{\tau \in \CW\,:\, |\tau| \ge 0\}$.
We interpret elements of $\hat \CT_-$
as triples $(F,\Labn,\Labe)$ as above, except that $F$ is now allowed to be a forest.
We also define the space $\CT_-$ by $\CT_- = \hat \CT_- / \CJ_+$.
With these notation at hand, we then define a map $\Deltam \colon \CT \to \CT_- \otimes \CT$ 
by setting, for $\tau = (T,\Labn,\Labe) \in \CW$,
\begin{equs}\label{e:Deltabar}
\Deltam \tau &= 
 \sum_{A \subset T} \sum_{\Labe_A,\Labn_A} \frac1{\Labe_A!}
\binom{\Labn}{\Labn_A}
 (A,\Labn_A+\pi\Labe_A, \Labe \restr E_A) 
 \\& \qquad  \otimes(\CR_A F,\CR_A (\Labn - \Labn_A), \Labe + \Labe_A)\;.
\end{equs}
Here, the sum runs over all $\Labn_A\colon V_A \to \N^2$ and 
$\Labe_A\colon \d(A,F) \to \N^2$, where $\d(A,F)$ denotes the edges in 
$E_F \setminus E_A$ that are adjacent to $V_A$. Also, for a function $\m \colon S \to \Z^2$
with $S$ a finite set, we write $\m! = \prod_{x \in S}\m(x)_0!\m(x)_1!$ and similarly
for the binomial coefficients, with the convention that $k!=\infty$ for $k < 0$. 
As before, the sum appearing here is actually finite
because all but finitely many terms have the first factor vanish in $\CT_-$.

Our motivation for the definition of $\Deltam$ is as follows. Assigning a number
to each $\tau \in \CW$ with $|\tau|<0$ is equivalent to choosing an algebra
morphism $g \colon \CT_- \to \R$. If we ignore for a moment the labels $\Labn$ and $\Labe$,
an operation of the type $M_g \colon \CT \to \CT$ with 
\begin{equ}[e:defMg]
M_g \tau = (g \otimes \id)\Deltam \tau\;,
\end{equ}
then corresponds to iterating over all ways of contracting subtrees of negative degree
contained in $\tau$ and replacing them by the corresponding constant assigned to it by $g$.
This corresponds to replacing a kernel of possibly several variables by a multiple
of a Dirac delta function forcing all arguments to collapse. The seemingly complicated
combinatorics appearing in \eqref{e:Deltabar} then encodes the possibility to 
also replace it by higher order derivatives of such a delta function in all of its arguments.

Similarly to before, $\Deltam$ can also be viewed as a map
$\Deltam \colon \CT_- \to \CT_- \otimes \CT_-$ by extending it multiplicatively 
from $\CW$ to all of $\CT_-$. 
Writing $\bullet$ for the product in $\CT_-$ (which has nothing to do with the product 
we have on $\CW$!), the following analogue to Proposition~\ref{prop:Hopf+}
was shown in \cite{YvainLorenzo}.

\begin{proposition}\label{prop:Hopf-}
There exists an algebra morphism $\CA_- \colon \CT_- \to \CT_-$ so that 
$(\CT_-, \bullet,\br\Deltam,\CA_-)$ is a Hopf algebra, and 
the map $\Deltam$, viewed as a map $\CT \to \CT_- \otimes \CT$, turns
$\CT$ into a left comodule for $\CT_-$.
\end{proposition}

If we write just as before $\CG_-$ for the group of characters 
$g \colon \CT_- \to \R$, this yields a right action of $\CG_-$ onto
$\CT$ by $g \mapsto M_g$ with $M_g$ as in \eqref{e:defMg}.
The following is then the main result of \cite{YvainLorenzo}.

\begin{theorem}\label{theo:algebra}
For every $g \in \CG_-$, the map $M_g$ is an admissible renormalisation procedure.
\end{theorem}

The idea of the proof of this theorem goes as follows. Assume for a moment
that one can also find a map $\Deltam \colon \hat \CT_+\to \hat \CT_- \otimes \hat \CT_+$
such that  
\begin{equ}
\Deltam \CJ_- \subset \CJ_+ \otimes \hat \CT_+ + \hat \CT_- \otimes \CJ_-\;,
\end{equ}
with $\CJ_\pm$ defined before \eqref{e:deltap1} and \eqref{e:Deltabar}. 
In particular, $\Deltam$ passes through
the quotients to a map $\CT_+ \to \CT_- \otimes \CT_+$, which we assume to satisfy the following.
\begin{claim}
\item On $\CT$, one has the identity
\begin{equ}[e:propWanted1]
\CM_-(\Deltam \otimes \Deltam) \Deltap = (\id \otimes \Deltap) \Deltam\;,
\end{equ}
where 
\begin{equ}[e:defM-]
\CM_- \colon \CT_-\otimes \CT \otimes \CT_- \otimes \CT_+
\to \CT_- \otimes \CT \otimes \CT_+
\end{equ}
is the map that multiplies the two factors in $\CT_-$.
The same is also true on $\CT_+$.
\item On $\CT_+$, one has the identity
\begin{equ}[e:propWanted2]
\Deltam \hat \CA_+ = (\id \otimes \hat \CA_+)\Deltam\;.
\end{equ}
\item The actions of $\CG_-$ onto $\CT$ and $\CT_+$ given by \eqref{e:defMg}
and the analogous formula for $\CT_+$ only increase degrees.
\end{claim}
In this case, it is straightforward to verify that, for any $g \in \CG_-$, if
we write $M_g$ as before, set $\PPi^g = \PPi M_g$ for some $\PPi$ such that
$\CZ(\PPi) = (\Pi,\Gamma)$ is a model, write $\CZ(\PPi^g) = (\Pi^g, \Gamma^g)$, and define $\gamma_{z\bar z}$
and $\gamma^g_{z\bar z}$ as in \eqref{e:defModel2}, one has
\begin{equ}
\gamma_{z\bar z}^g = (g \otimes \gamma_{z\bar z})\Deltam\;,\qquad 
\Pi_z^g  = (g \otimes \Pi_z)\Deltam\;.
\end{equ}
To show this, one first uses \eqref{e:propWanted2} to show that 
$f_{z}^g = (g \otimes f_z)\Deltam$,
where $f$ and $f^g$ are defined from $\PPi$ and $\PPi^g$ as in \eqref{e:defModel1}.
One then uses \eqref{e:propWanted1} (on $\CT$) to show that the required identity for $\Pi_z^g$ holds. 
Finally, one uses \eqref{e:propWanted1} on $\CT_+$ to show that if 
one views $M_g$ as acting on
$\CT_+$ via \eqref{e:defMg}, then its action distributes over the product defined
in \eqref{e:defProd} in the sense that $(M_g f)\circ (M_g \bar f) = M_g (f\circ \bar f)$,
which then implies the required identity for $\gamma_{z\bar z}^g$. The fact that
the action of $M_g$ increases degrees guarantees that 
$\CZ(\PPi^g)$ is again a model, provided that $\CZ(\PPi)$ is.

The problem is that \eqref{e:propWanted1} actually fails in our situation. 
However, it turns out that it can still be rescued by the following construction.
We look for a larger space $\CT^\ex$, together with corresponding spaces
$\hat \CT^\ex_+$, $\hat \CT^\ex_-$ and ideals $\CJ_-^\ex \subset \hat \CT^\ex_+$
and $\CJ_+^\ex \subset \hat \CT^\ex_-$,
all of them $\R$-graded, 
such that the following properties hold.
\begin{claim}
\item There are analogous maps to $\Deltap$ and $\Deltam$ acting on these ``extended''
spaces and
such that all of the algebraic relations described above are satisfied, including
\eqref{e:propWanted1} and \eqref{e:propWanted2}. In particular, one has a map 
$\CZ^\ex$ turning linear maps $\CT^\ex \to \CC^\infty$ into candidate models 
on the regularity structure $(\CT^\ex,\CG_+^\ex)$ defined in formally the same
way as above.
\item There exists a projection $\pi^\ex\colon \CT^\ex \to \CT$ which is a right inverse
for the inclusion $\CT \hookrightarrow \CT^\ex$ and is such that 
for any admissible $\PPi\colon\CT \to \CC^\infty$, 
$\CZ(\PPi)$ is a model if and only if $\CZ^\ex(\PPi \pi^\ex)$ is a model.
\item There exists an algebra morphism $\pi^\ex_-\colon \CT_-^\ex \to \CT_-$
such that, for every $g \in \CG_-$, one has 
\begin{equ}[e:proppi]
\pi^\ex M^\ex_{g\pi^\ex_-} = M_g \pi^\ex\;.
\end{equ}
\end{claim} 
Once we have constructed these larger spaces, the proof of Theorem~\ref{theo:algebra} is
rather straightforward. Fix $g\in \CG_-$ and $\PPi$ such that $\CZ(\PPi)$ is
a model. Then, by the second property above, in order to show that $\CZ(\PPi M_g)$
is a model, it suffices to show that $\CZ^\ex(\PPi M_g \pi^\ex)$ is a model.
However, by \eqref{e:proppi}, we have $\CZ^\ex(\PPi M_g \pi^\ex) = \CZ^\ex(\PPi \pi^\ex M^\ex_{g\pi^\ex_-})$
and, again by the second property, we already know that $\CZ^\ex(\PPi \pi^\ex)$ is
a model. We conclude by the fact that $M^\ex_{g\pi^\ex_-}$ is an admissible renormalisation
map thanks to the argument given above, using the properties \eqref{e:propWanted1} 
and \eqref{e:propWanted2} for the maps $\Deltapm$ defined on the extended spaces.

\section{Analytical aspects}
\label{sec:analysis}

At this stage, we have built a rather large group $\CG_-$ acting on 
our space of formal expressions $\CT$ by admissible renormalisation
procedures. Consider now regularised space-time noises $\xi_i^{(\eps)}$
as in \eqref{e:scalerho} and define $\PPi^{(\eps)}$ as their canonical
lift, defined by \eqref{e:model}. The following result is a particular
instance of the main theorem of \cite{Ajay}.

\begin{theorem}\label{theo:BPHZ}
There exists a choice of (deterministic) elements $g^\eps \in \CG_-$ such that,
setting $M^\eps = M_{g^\eps}$ as in \eqref{e:defMg}, the sequence of models
$\CZ(\PPi^{(\eps)} M^\eps)$ converges to a limiting model in $\MM$.
\end{theorem}

Before we give an idea of the proof of this theorem, let us show how the
element $g^\eps$ determining the suitable renormalisation procedure is constructed.
It turns out that this construction is very similar to the construction of the 
elements $f_z \in \CG_+$ constructed in \eqref{e:defModel1} in order to correctly recenter 
the model, so that its behaviour around a given space-time point $z$ matches
its degree. This is maybe not surprising since one can also view the renormalisation as some
kind of ``recentering procedure'' except that this time, instead of insisting that the
\textit{evaluation} of the model at a given location vanishes for basis vectors of
positive homogeneity, we would like to impose that the \textit{expectation} of the model
vanishes for basis vectors of negative homogeneity.

Recall equation \eqref{e:Deltabar} defining $\Deltam$. It follows from this definition
that the antipode $\CA_-$ for the Hopf algebra $\CT_-$ is defined 
for $\tau = (F,\Labn,\Labe)$ by the recursion
\begin{equ}
\CA_- \tau = -
 \sum_{A \subset T \atop A \neq \emptyset} \sum_{\Labe_A,\Labn_A} \frac1{\Labe_A!}
\binom{\Labn}{\Labn_A}
 (A,\Labn_A+\pi\Labe_A, \Labe \restr E_A) 
   \bullet \CA_- (\CR_A F,\CR_A (\Labn - \Labn_A), \Labe + \Labe_A)\;,
\end{equ}
where $\bullet$ is the product in $\CT_-$ as before. As above, the renormalisation
procedure involves a twisted antipode. In order to define this, recall that 
$\CT_- = \hat \CT_- / \CJ_+$ and write
$\pi_- \colon \hat \CT_- \to \hat \CT_-$ for the projection onto elements of strictly
negative degree. Similarly to above, we then define $\hat \CA_- \colon \CT_- \to \hat \CT_-$ 
inductively as being the unique algebra morphism so that, 
on elements of the type $(F,\Labn,\Labe)$ with $F$ a single tree (so that it is 
identified with an element of $\CW$), one has
\begin{equs}
\hat \CA_- (F,\Labn,\Labe) &= 
 -\sum_{A \subset T \atop A \neq \emptyset} \sum_{\Labe_A,\Labn_A} \frac1{\Labe_A!}
\binom{\Labn}{\Labn_A} \\
& \pi_- \bigl((A,\Labn_A+\pi\Labe_A, \Labe \restr E_A) 
   \bullet \hat \CA_- (\CR_A F,\CR_A (\Labn - \Labn_A), \Labe + \Labe_A)\bigr)\;.
\end{equs}
We furthermore note that any \textit{random} linear map
$\PPi\colon \CT \to \CC^\infty$ with finite expectation gives rise to a character
$g^-(\PPi)$ on $\hat \CT_-$ by simply setting
\begin{equ}
g^-(\PPi)\tau = \E (\PPi \tau)(0)\;,\qquad \tau \in \CW\;,
\end{equ}
and then extending it multiplicatively.
In this setting, at least for instances of $\PPi$ satisfying a suitable kind 
of stationarity, we then claim that the ``correct'' choice of character $g^\eps$
appearing in Theorem~\ref{theo:BPHZ} is given by
\begin{equ}
g^\eps  = g^-(\PPi^{(\eps)})\hat \CA_-\;.
\end{equ}
Comparing this and \eqref{e:defMg} to \eqref{e:defModel1}, we see that
the renormalisation procedure required to make our models converge to a 
finite limit is indeed formally identical (modulo changing $+$ / evaluation 
into $-$ / expectation)
to the recentering procedure discussed before.

Combining these constructions with the results of \cite[Sec.~10]{Regularity},
one concludes that Theorem~\ref{theo:BPHZ} essentially 
follows as soon as one has an estimate of the type
\begin{equ}[e:estimate]
\E \scal{\bigl(g^-(\PPi^{(\eps)})\hat \CA_- \otimes \PPi^{(\eps)}
\otimes g_z^+(\PPi^{(\eps)})\hat \CA_+\bigr)\CM_-(\Deltam\otimes \Deltam)\Deltap\tau, \phi_z^\lambda}^2 \lesssim \lambda^{2|\tau|}\;,
\end{equ}
for every $\tau \in \CW$ with $|\tau| < 0$, where $|\tau|$ denotes the degree of
$\tau$ as before, $\CM_-$ is as in \eqref{e:defM-}, 
$\scal{\cdot,\cdot}$ is the usual scalar product in  $L^2(\R^2)$,
and $\phi_z^\lambda$ is the translate and rescaling of a test function $\phi \in \CC_0^\infty$
with sufficiently many derivatives bounded by $1$ and with support in the ball of radius $1$
around the origin as in \cite[Def.~3.3]{Review}.

From an algebraic perspective, the definitions of $\Deltam$ and $\hat \CA_-$ are very
strongly reminiscent of the Hopf-algebraic formulation of Zimmermann's forest formula
\cite{Kreimer}, which was further explored in \cite{CK,CK2}.
More precisely, our space $\hat \CT_-$ is analogous to the space $\CA$ in \cite{Kreimer},
the quotiented space $\CT_-$ is analogous to the space $\CA/\sim$, etc.,
so \eqref{e:estimate} is really a form of BPHZ theorem. 

The difference is threefold.
First, our basic combinatorial structure is described by collections of trees rather than
Feynman diagrams. These can then be interpreted as generating 
Feynman diagrams when taking expectations, by contracting
leaves according to Wick's theorem. (Or the cumulant formula if one considers a $\PPi^{(\eps)}$
which is not obtained from the lift of a Gaussian process.)
Second, the result in \cite{Ajay} applies to very large class of kernels $K$, provided that they
exhibit the ``correct'' behaviour near the origin and, unlike most related results that
can be found in the literature, it doesn't rely on the driving noise being Gaussian. 
Finally, and this is really the main difference,
we see both ``positive'' and ``negative'' renormalisations appearing in
\eqref{e:estimate}, while the usual calculations performed in the context of QFT only
involve negative renormalisation. The purpose of the latter is to ensure that we obtain
finite quantities in the limit $\eps \to 0$. The former on the other hand is crucial in order 
to obtain the correct power of $\lambda$ in the right hand side of \eqref{e:estimate}.

Let us explain the main ingredients appearing in the proof of \eqref{e:estimate}.
First, by translation invariance, one can set $z=0$.
Using Wick's formula, the left hand side of \eqref{e:estimate} can then be written,
for some $N > 2$, as 
\begin{equ}
I_\tau^{\lambda,\eps} = \int \phi^\lambda(z_1)\phi^\lambda(z_2) \,\CK_\tau^{(\eps)}(z_1,\ldots,z_N)\,dz_1\cdots dz_N\;,
\qquad z_i \in \R^2\;,
\end{equ}
for some smooth kernel $\CK_\tau^{(\eps)}$ which, as $\eps \to 0$, converges 
to a smooth limit $\CK_\tau$ on the configuration space $C_N = \{(z_1,\ldots,z_N)\,:\, z_i \neq z_j\,\forall\,i\neq j\}$, but exhibits quite singular behaviour on the ``big diagonal'' where
two or more arguments coincide.

In order to estimate an integral of this type, we use the following construction
reminiscent of the Fulton-MacPherson compactification of $C_N$ and already
used in \cite{HS,HQ}.
First, note that if $\Lambda_\lambda \in (\R^2)^N$ denotes the support of 
$z\mapsto \phi^\lambda(z_1)\phi^\lambda(z_2)$, $\Trees$ is a countable index set, and
$\{D_\sigma\}_{\sigma \in \Trees}$ are a collection of bounded regions exhausting
all of $(\R^2)^N$, an integral as above can trivially be estimated by
\begin{equ}[e:boundI]
|I_\tau^{\lambda,\eps}| \lesssim \lambda^{-6} \sum_{\sigma \in \Trees} \one_{D_\sigma \cap \Lambda_\lambda \neq \emptyset} |D_\sigma|
\sup_{z \in D_\sigma} |\CK_\tau^{(\eps)}(z)|\;,
\end{equ}
where $|D|$ denotes the Lebesgue measure of $D$ and the factor $\lambda^{-6}$ comes
from the parabolic rescaling of $\phi^\lambda$.
Such a bound is close to optimal if the regions $D_\sigma$ can be chosen in
such a way that the integrand $\CK_\tau^{(\eps)}$ does not vary much over them.

A good index set $\Trees$ turns out to be given by the set of rooted binary trees $T = (V_T,E_T)$
with $N+1$ leaves endowed with an integer label $n(x)$ at each interior vertex $x \in V_T$.
We now build a map $\CS\colon C_N \to \Trees$ such that the sets $D_\sigma$ are then
given by $D_\sigma = \CS^{-1}(\sigma)$. 
Given $z \in C_N$, $\CS(z)$ is built by simultaneously looking at a sequence $\CP_k$ of
partitions of $[N] = \{0,\ldots,N\}$ and $(T_k,n_k) = (V_k, E_k, n_k)$ of labelled 
graphs (with labels in $\Z \cup \{+\infty\}$) as follows. 
We set $\CP_0 = \{\{0\},\ldots,\{N\}\}$, $V_0 = [N]$, $E_0 = \emptyset$, and $n_0(x) = +\infty$ 
for $x \in V_0$. Then, given $\CP_j$ and $(T_j,n_j)$, we define the next element in the
sequence as follows. If $\CP_j = \{[N]\}$, then the construction stops
and we set $\CS(z) = (T_j,n_j)$. Otherwise, for every $A \in \CP_j$, consider the set
$z_A = \{z_i\}_{i \in A}$ (where we set $z_0 = 0$) and consider the pair $(A,B)$ such that $d(z_A,z_B) \le d(z_C,z_D)$
for every $C,D \in \CP_j$, where $d$ denotes the Hausdorff distance between compact sets. 
Since the points $z_i$ are all distinct, the pair $(A,B)$ is unique.
One then sets
\begin{equs}
\CP_{j+1} &= (\CP_j \setminus \{A,B\}) \cup \{A\cup B\}\;,\\ 
V_{j+1} &= V_j \cup \{A\cup B\}\;, \quad
E_{j+1} = E_j \cup \{(A\cup B,A), (A\cup B,B)\}\;.
\end{equs} 
We furthermore define $n_{j+1}$ to be equal to $n_j$ on $V_j$
and $n_{j+1}(A\cup B)$ to be the only integer such that 
\begin{equ}
d(z_A,z_B) \le 2^{-n_{j+1}(A\cup B)} < 2d(z_A,z_B)\;.
\end{equ}
As a consequence of the properties of the Hausdorff distance, the labelled trees
produced in this way have the property that if we partially order $T$ 
in the natural way so that the root $[N]$
is minimal, it is always the case that  $x \le y$ implies $n(x) \le n(y)$.
The following lemma is crucial.

\begin{lemma}\label{lem:props}
For $z \in C_N$, let $(T,n) = \CS(z)$. Then, there exist constants $c$ and $C$
depending only on $N$ such that, for any $i,j \in [N]$,
and writing $i \wedge j = \sup\{x\,:\, x\le \{i\} \;\&\; x \le \{j\}\}$ 
for the most recent common ancestor of 
the leaves $\{i\}$ and $\{j\}$ in $T$, one has
\begin{equ}[e:distances]
c2^{-n_{j+1}(i\wedge j)} \le d(z_i,z_j) \le C2^{-n_{j+1}(i\wedge j)}\;.
\end{equ}
(Again with the convention $z_0 = 0$.) Furthermore, for $\sigma = (V,E,n)$, 
one has the upper bound
\begin{equ}[e:boundD]
|D_\sigma|
\lesssim \prod_{x \in V_\circ} 2^{-d n(x)}\;,
\end{equ}
uniformly over $\sigma \in \Trees$, where $V_\circ \subset V$ denotes the set of those
vertices that are not leaves in $T$.
\end{lemma}

In other words, modulo constant factors, the distance between any two points of
a given configuration $z$ is completely determined by $\CS(z)$. Furthermore, the tree
structure provides a very efficient way of encoding the various
constraints given by the triangle inequality.
The aim then is to obtain, for any such binary tree $T = (V,E)$
a function $\bar\eta_T \colon V \to \R$, such that the integrand $\CK_\tau^{(\eps)}$ can be bounded by
\begin{equ}[e:wantedK]
|\CK_\tau^{(\eps)}(z)| \lesssim \prod_{x \in V_\circ} 2^{-\bar \eta_T(x) n(x)}\;,
\end{equ}
uniformly in $\eps$. Combining this with \eqref{e:boundD} and \eqref{e:distances},
one then obtains a bound of the type
\begin{equ}[e:boundI]
|I_\tau^{\lambda,\eps}| \lesssim \lambda^{-6} 
\sum_{T=(V,E)}\,\sum_n \prod_{x \in V_\circ}2^{-\eta_T(x)\, n(x)}\;,
\qquad \eta_T(x) = \bar \eta_T(x) + d\;,
\end{equ}
where the inner sum runs over all weakly increasing maps $n \colon V_\circ \to \Z$
such that furthermore 
\begin{equ}[e:condn]
2^{-n(0\wedge 1\wedge 2)} \le C \lambda\;,
\end{equ}
for some fixed $C > 1$. (Here $i\wedge j$ is as in Lemma~\ref{lem:props}.) 
This encodes the constraint that one only considers terms
in the sum such that $D_\sigma \cap \Lambda_\lambda \neq \emptyset$.
It is relatively straightforward to verify recursively that one has indeed a bound
of the type \eqref{e:wantedK} so that \eqref{e:boundI} holds with functions 
$\eta_T$ satisfying $\sum_{x\in V_\circ} \eta_T(x) = |\tau|+6$.
It therefore remains to obtain conditions on such functions $\eta_T$ guaranteeing that
a sum of the type appearing in \eqref{e:boundI} is bounded by 
$C \lambda^{\sum_{x\in V_\circ} \eta_T(x)}$ for some fixed $C$.

This is the content of the following lemma.

\begin{lemma}\label{lem:sum}
Let $T = (V,E)$ be a rooted binary tree with leaves equal to $[N]$ with $N > 2$ 
and let $\eta \colon V_\circ \to \R$. 
Set $x_\star = 0\wedge 1 \wedge 2 \in V_\circ$ and let $V_\star$
denote the nodes $x \in V_\circ$ lying on the path from $x_\star$ to the root, 
but not the root itself. Assume furthermore that the following two conditions hold.
\begin{claim}
\item[1.] For every $x \in V_\circ$ one has
$\sum_{y \ge x} \eta(y) > 0$.
\item[2.] For every $x \in V_\star$ one has $\sum_{y \not \ge x} \eta(y) < 0$.
\end{claim}
Then, there exists a constant $C$ such that 
\begin{equ}
\sum_n \prod_{x \in V_\circ}2^{\eta(x)\, n(x)} \le C \lambda^{-\sum_{x\in V_\circ} \eta(x)}\;,
\end{equ}
where the sum runs over all increasing $n \colon V_\circ \to \Z$ satisfying \eqref{e:condn}.
\end{lemma}

The proof of this lemma is relatively straightforward and can be found
for example in \cite{HQ}. The two conditions appearing here are generalisations
of the standard conditions on the integrability of a function with power law behaviour
at the origin (corresponding to the first condition) and at infinity
(corresponding to the second condition). The problem in our case is that if we simply 
replaced the complicated expression appearing in the left hand side of \eqref{e:estimate}
by $\E \scal{\PPi^{(\eps)}\tau, \phi_z^\lambda}^2$, then although 
\eqref{e:wantedK} and \eqref{e:boundI} would still be satisfied for
some $\eta_T$ with $\sum_{x\in V_\circ} \eta_T(x) = |\tau|+6$, both  
conditions of Lemma~\ref{lem:sum} would fail for a typical $\tau \in \CW$.

The purpose of the two renormalisation procedures encoded by $\Deltam$ and $\Deltap$ 
appearing in \eqref{e:estimate} is to allow us to obtain
an improved bound which involves a function $\eta_T$ that does satisfy the conditions
of Lemma~\ref{lem:sum}.
In this procedure, the purpose of the ``negative renormalisation'' is precisely to 
ensure that the first condition is satisfied (thus removing small-scale divergences),
while the purpose of the ``positive renormalisation'' is to ensure that 
the second condition is satisfied, guaranteeing integrability at large scales.
The combinatorics of overlapping divergencies can in particular  be unraveled 
by adapting tools from \cite{FMRS}.

\subsection*{Acknowledgements}

{\small
It is a pleasure to thank Yvain Bruned, Ajay Chandra and Lorenzo Zambotti for
our numerous discussions and debates in the process of obtaining these results.
This research was supported by the Leverhulme Trust through a leadership award 
and by the ERC through a consolidator award.
}

\bibliographystyle{Martin}

\bibliography{proceedings}

\end{document}